\documentclass[12pt]{article}
\usepackage{amsmath,amsthm,amsfonts,amssymb,amscd}
\usepackage{graphics}
\usepackage[utf8]{inputenc}
\usepackage{graphicx}
\usepackage{color}
\newcommand {\dem}{{\noindent\it Proof:}}
\newcommand {\fimdem}{\rightline{$\square$}}
\newcommand{\w}[1]{\widehat{#1}}

\newcommand{\dis}{\displaystyle}
\newcommand {\R}{\mathbb{R}}
\renewcommand{\over}[1]{\overline{#1}}

\renewcommand{\thefootnote}

\numberwithin{equation}{section}
\newtheorem{lema}{Lemma}[section]
\newtheorem{teo}{Theorem}[section]
\newtheorem{obs}{Remark}[section]
\newtheorem{prop}{Proposition}[section]

\begin{document}

\centerline{\Large\bf Asymptotic properties for second-order linear
evolution } \centerline{\Large\bf problems with fractional laplacian
operators}

\vskip .4in

\centerline{\sc Ma\'ira Fernandes Gauer Palma}

\centerline{\sl Department of Mathematics, }

\centerline{\sl Federal University of Santa Catarina}

\centerline{\sl Campus Universit\'ario, Trindade, 88040-900,
Florian\'opolis, SC, Brazil,}

\centerline{\sl  e-mail: maira.gauer@ufsc.br}

\vskip .2in

\centerline{\sc Cleverson Roberto da Luz}

\centerline{\sl Department of Mathematics, }

\centerline{\sl Federal University of Santa Catarina}

\centerline{\sl Campus Universit\'ario, Trindade, 88040-900,
Florian\'opolis, SC, Brazil,}

\centerline{\sl  e-mail: cleverson.luz@ufsc.br}

\footnote{\hspace{-0,75cm} 2010 Mathematics Subject Classification.
Primary: 35B40; Secondary:  35L15, 35K90, 35L90.
\newline Key words: Asymptotic behavior, fractional Laplace operator,
Fourier space, second-order equations.
\newline * The second author is supported by CNPq (Brazil), Proc. 308868/2015-3.} \vskip .3in

\begin{minipage}{14cm}
{\bf Abstract:} In this work we study the asymptotic behavior of
solutions for a general linear second-order evolution differential
equation in time with fractional Laplace operators in $\R^n$. We
obtain improved decay estimates with less demand on the initial data
when compared to previous results in the literature. In certain
cases, we observe that the dissipative structure of the equation is
of regularity-loss type. Due to that special structure, to get decay
estimates in high frequency region in the Fourier space it is
necessary to impose additional regularity on the initial data to
obtain the same decay estimates as in low frequency region. The
results obtained in this work can be applied to several initial
value problems associated to second-order equations, as for example,
wave equation, plate equation, IBq, among others.
\end{minipage}

\section{Introduction}
\hspace{1cm}We consider the following Cauchy problem with fractional
Laplace operators in $\R^n$:
\begin{equation}\label{eq}v_{tt}(t,x)+(-\Delta)^\delta v_{tt}(t,x)+(-\Delta)^\alpha v(t,x)
+ (-\Delta)^{\theta}v_t (t,x)= 0, \quad t\geq 0,\,\,\,\, x\in \R^n
\end{equation}
with initial data
\begin{equation}\label{dados}
v(0,x) = v_0(x), \quad v_t(0,x)=v_1(x),
\end{equation}
where $\delta, \alpha$ and $\theta$ are real numbers with 
$\delta, \alpha \geq 0$ and $\theta \in\left[0,\alpha\right]$. The fractional power operator
$\,(-\Delta)^\theta : H^{2\theta}({\R}^n) \subset L^2({\R}^n)
\rightarrow L^2({\R}^n)$ $(\theta \geq 0)\,$ is defined by
$$(-\Delta)^\theta v(x):= {\cal F}^{-1} \bigl(|\xi|^{2\theta} \widehat{v}(\xi)\bigl)(x),
\qquad v \in H^{2\theta}({\R}^n), \,\,\,\, x \in {\R}^n,$$ where
${\cal F}$ denotes the usual Fourier transform in $L^2({\R}^n)$ with
respect to the $x$ variable, $\widehat{v} = {\cal F}(v)$,
$H^s=H^s(\R^n)$ denotes the usual Sobolev space of $L^2$ functions
equipped with the norm $\|\,\cdot\,\|_{H^s}$ and $|\cdot|$ denotes
the usual norm in ${\R}^n$. The operator $(-\Delta)^\theta$ is
nonnegative and self-adjoint in $L^2({\R}^n)$.

The total energy $E_v(t)$ associated to the solution $v(t)$ of
equation (\ref{eq}) is defined by
\begin{equation}\label{energia}
E_v(t)=\frac{1}{2}\left\{ ||v_t(t)||^2+||(-\Delta)^{\frac{\delta}{2}} v_t(t)||^2
+||(-\Delta)^{\frac{\alpha}{2}}v(t)||^2 \right\}.
\end{equation}
For simplicity of
notations, in particular, we use in all text $\| \cdot \|$ instead
of $\| \cdot \|_{L^2}$.

In order to obtain estimates for the problem
(\ref{eq})-(\ref{dados}) we work with the \linebreak corresponding
problem in the Fourier space. Applying the Fourier transform with
\linebreak respect to the variable $x$, we obtain
\begin{equation}\label{f}
\left\{ \begin{array}{ll}
(1+|\xi|^{2\delta})\,\widehat{v}_{tt}(t,\xi)+
|\xi|^{2\theta}\,\widehat{v}_t(t,\xi)+|\xi|^{2\alpha}\,\widehat{v}(t,\xi)=0,
\quad t\geq 0,\,\,\,\, \xi\in \R^n
\\
\widehat{v}(0,\xi)=\widehat{v}_0(\xi), \quad
\widehat{v}_t(0,\xi)=\widehat{v}_1(\xi), \quad \xi \in \R^n.
\end{array} \right.
\end{equation}

The eigenvalues of the problem (\ref{f}) have nonpositive real part,
and they are given by
$$\lambda_{\pm}=\frac{|\xi|^{2\theta}}{2(1+|\xi|^{2\delta})}\left(-1
\pm \sqrt{1-4|\xi|^{2(\alpha-2\theta)}(1+|\xi|^{2\delta})}\right).$$

The solution of (\ref{f}) can be written as
$$\w{v}(t,\xi)=\w{K}_0(t,\xi) \w{v}_0(\xi) + \w{K}_1(t,\xi)
\w{v}_1(\xi),$$
where
\begin{equation}\label{K0K1}
\w{K}_0(t,\xi)=\dfrac{\lambda_{+}e^{\lambda_{-}t}
-\lambda_{-}e^{\lambda_{+}t}}{\lambda_{+}-\lambda_{-}}
\hspace{0,8cm} \text{and} \hspace{0,8cm}
\w{K}_1(t,\xi)=\dfrac{e^{\lambda_{+}t}-e^{\lambda_{-}t}}{\lambda_{+}-\lambda_{-}}.
\end{equation}
Therefore
$$v(t,x)=kK_0(t,x)\ast v_0(x) +k K_1(t,x)\ast v_1(x),$$
with $k$ a positive constant such that $\w{f} \cdot \w{g}=k(f\ast
g)$.

For $\varepsilon>0$, we denote by $E_0(t,x)$ and $E_\infty(t,x)$ the
solution of (\ref{eq}) localized to low and high frequencies, that
is,
\begin{eqnarray}
& E_0=E_0(t,x)(v)=\mathcal{F}^{-1}(\chi(\xi)\w{v}(t,\xi)),
\label{Ezero}
\\
&
E_\infty=E_\infty(t,x)(v)=\mathcal{F}^{-1}((1-\chi(\xi))\w{v}(t,\xi)),\label{Einf}
\end{eqnarray}
where $\chi(\xi)$ is the characteristic function of $\{ \xi \in
\R^n\,\,/\,\,\, |\xi| < \varepsilon \}$. To estimate the \linebreak
$x$-derivatives norms of $v$ and $v_t$ it is sufficient to estimate
the derivatives norms of $E_0$, $\partial_tE_0$, $E_\infty$ and
$\partial_tE_\infty$.

If $\theta < \delta$, we observe that the decay structure of
(\ref{eq}) is of regularity-loss type which is characterized by the
structure of the eigenvalues associated to the problem. The
regularity-loss property ceases to occur in the case when $\theta =
\delta$. Due to that special structure, to get decay estimates in
the high frequency region in the Fourier space it is necessary to
impose additional regularity on the initial data to obtain the same
decay estimates as in the low frequency region. If $\delta \leq
\theta$ this effect does not appear, since the solution decays
exponentially in the high frequency zone of the Fourier space (see
Proposition \ref{alta}). Such decay property of the regularity-loss
type was also investigated for the dissipative Timoshenko system
\cite{KISK}, the plate equation under rotational inertia
effects in $\R^n$ \cite{SK, RCR2, DRC1} and a hyperbolic-elliptic
system of a radiating gas model \cite{TKSK}.

Estimates for the solution to the wave equation with structural
damping
\begin{equation}\label{we}
u_{tt}(t,x) -\Delta u(t,x)  + 2a \,(-\Delta)^\theta u_t(t,x) =
0\,,\qquad t\geq0, \,\,\quad x\in\R^n,
\end{equation}
have been derived in \cite{IN, RCR1, DR, DE-2,karch}, where as
the case of time-dependent damping coefficients is considered in
\cite{Wi, LR}.  In \cite{DE-2}, the authors decompose the
solution to (\ref{we}) into two parts, $u= u^+ + u^-$, each one
related to one of the two characteristics roots associated to
(\ref{we}). The asymptotic behavior of the Fourier transforms of
each part hints to two different diffusion phenomena. This type of
problem have been extensively studied in the mathematical and
physical literature (see, for instance, \cite{BKW, CC, FK, JV}).
Karch \cite{karch} studied the large time behavior of solutions to
the initial value problem (\ref{eq}) with $\delta=0$, $0 \leq 2
\theta < \alpha$ and a non-linear term $F(x,t,u,u_t,\nabla u)$. In
the cited paper, an analysis of the solution formula of the linear
problem leads to the conclusion that they behave, as $t \rightarrow
\infty$, like solutions of a similar diffusion equation to the
problem in \cite{DE-2}. 

Next we mention some important previous works which are related with
this paper. Sugitani-Kawashima \cite{SK} studied a semilinear
dissipative plate equation with rotational inertia effects and a
frictional dissipation. To the linear problem, they used the
explicit solution and equivalences for eigenvalues. To the
semilinear problem, they introduce a Banach space $X$ defined by the
norm
$$||u||_X=\sum_{\sigma_0(k)\leq s+1} \sup_{t\geq 0} (1+t)^{\frac{k}{4}}||\partial_x^k
u(t)||_{H^{s+1-\sigma_0(k)}} $$ where $\sigma_0(k)=k+\left[\frac{k+1}{2}\right]$
and $s\geq \sigma_0(k)-1$.
They proved that the problem is solved globally in time in the above
function space and found optimal decay estimates of solutions under
the additional regularity assumption on the initial data $u_0\in
H^{s+1}$ and $u_1\in H^s$, for $s$ large enough.

D'Abbicco-Reissig \cite{DR} studied the semilinear structural damped
wave equation. They obtained optimal decay to the norm of solution,
to the energy terms and determined the influence of fractionary
dissipation on the critical exponent. The method used for the linear
problem is similar to the method applied in \cite{SK} (see also
\cite{Wi}).

Based on the energy method in Fourier space, introduced by
Umeda-Kawashima-Shizuta \cite{UKS}, several works (see
\cite{DRK,IN,ITY}) showed decay estimates to some evolution
equations in $\R^n$. Char\~ao-da Luz-Ikehata \cite{RCR1, RCR2}
introduced a new method to get decay estimates supported on the
energy method in Fourier space combined with the Haraux-Komornik
inequality, the monotonicity of the energy density in the Fourier
space and the property of power singularities less than $n$ are
integrable around the origin of $\R^n$. They obtained in \cite{RCR1}
almost optimal decay for wave equation with a fractional damping and
in \cite{RCR2} for the plate equation with rotational inertia
effects and fractional damping. In \cite{RCR} they extended these
results to an abstract problem of second order differential
equation. The decay rate $E(t)= O(t^{-k})$ is almost optimal  means
that
$$E(t) = O(t^{-k +\varepsilon})
\qquad (t \rightarrow +\infty)$$ for any $\varepsilon
> 0$.

In this work we study a more general equation when compared to the
problems (or at least to the associated linear problem) studied in
\cite{RCR1, RCR2, IN, DE-2, DR, karch, SK, PE, wang2, wang}. The main objective of this work is to improve the results
obtained in \cite{RCR}. Due to the method used in \cite{RCR}, the
authors proved almost optimal decay rates only for the total energy
of order $\alpha$ of the problem, that is, the method does not allow
to obtain optimal decay rates for the total energy, as well as does
not allows to estimate each term of the energy separately. However,
as can be seen in Section 4 of this work, the rates can vary for
each term of the energy. In this paper we obtain optimal rates of
decay for $\|\partial_x^{\gamma_1}v(t)\|$,
$\|\partial_x^{\gamma_2}v_t(t)\|$ and the corresponding regularity
of the initial data.

To obtain the desired decay rates, we study low and high frequency
regions \linebreak separately. The next section will deal low
frequency, using the ideas presented in \cite{DR} and \cite{SK},
i.e., we consider the explicit solution of the problem and estimate
the \linebreak eigenvalues. It will be necessary to separate in two
cases: real eigenvalues ($\alpha >2\theta$) and complex eigenvalues
($\alpha \leq 2\theta$).

In the third section we study the problem in high frequency zone
through a redesign of the energy method in Fourier space introduced
by Char\~ao-da Luz-Ikehata in \cite{RCR1, RCR2}. In this section we
find the additional regularity required in the initial data (if $
\theta <\delta $) to obtain desired decay rates.

In the last section we exhibit and prove the main theorems,
combining the results of Sections 2 and 3. In addition, some
important applications are presented, finding optimal decay rates
for the norm of solution and for the terms of the energy associated
with the wave equation with fractional damping, plate equation with
rotational inertia effects and fractional damping and a Boussinesq
equation with fractional damping.

The method may further be applied to various other evolution
equations in $\R^n$ with constant coefficients. Moreover, we can add
terms of the type $(-\Delta)^{\delta_1}v_{tt}$,
$(-\Delta)^{\theta_1}v_t$, $(-\Delta)^{\alpha_1}v$ in equation
(\ref{eq}) and get decay rates to the solution of the new problem
from the results obtained for the problem (\ref{eq})-(\ref{dados}).
This is possible because, for example, adding the term
$(-\Delta)^{\delta_1}v_{tt}$, we obtain
$$|\xi|^{2\delta}+|\xi|^{2\delta_1}$$ as coefficient of the term
$\w{v}_{tt}$ at the equation in Fourier space. But this coefficient
is equivalent to $|\xi|^{2\delta_0}$, where $\delta_0=\min\{\delta,
\delta_1\}$ in low frequency, and $\delta_0=\max\{\delta,\delta_1\}$
in high frequency. The same occurs when we add terms of the type
$(-\Delta)^{\theta_1}v_t$, $(-\Delta)^{\alpha_1}v$ in equation
(\ref{eq}). This means that using the method presented in this work,
it is also possible to obtain estimates for this new problem. In
order not to make this a long work, we're going to consider only a
simple example in order to illustrate this case (see Subsection
4.2.3).

Note that in this work we do not discuss the existence and
uniqueness of solutions because they can be obtained via standard
methods.

%%%%%%%%%%%%%%%%%%%%%%%%%%%%%%%%%%%%%%%%%%%%%%%%%
%%%%%%%%%%%%%%%%%%%%%%%%%%%%%%%%%%%%%%%%%%%%%%%%%
%%%%%%%%%%%%%% BAIXA FREQUENCIA %%%%%%%%%%%%%%%%%
%%%%%%%%%%%%%%%%%%%%%%%%%%%%%%%%%%%%%%%%%%%%%%%%%
%%%%%%%%%%%%%%%%%%%%%%%%%%%%%%%%%%%%%%%%%%%%%%%%%

\section{Low frequency region: $|\xi|<\varepsilon$}

\hspace{1cm}In the following we use the notation $f\lesssim g$ to
mean that $0\leq f\leq Cg$ for some constant $C>0$. The notation $g
\approx f$ means $g \lesssim f$ and $f \lesssim g$. Furthermore, we
use the symbol to $ \gamma $ denote a multi-index with non-negative
entries.

We observe that
\begin{eqnarray}
\nonumber||\partial_t^j \partial_x^{\gamma} E_0(t)||^2 &&
\hspace{-0,5cm}=\int_{\R^n}|\partial_t^j
\w{v}(t)|^2|\xi|^{2|\gamma|} \chi(\xi)^2 d\xi\\\nonumber &&
\hspace{-0,5cm}\lesssim ||\w{v}_0||_{L^\infty}^2\int_{\R^n}
|\partial_t^j
\w{K}_0(t)|^2|\xi|^{2|\gamma|}\chi(\xi)^2d\xi+||\w{v}_1||_{L^\infty}^2\int_{\R^n}
|\partial_t^j \w{K}_1(t)|^2|\xi|^{2|\gamma|}\chi(\xi)^2d\xi,
\end{eqnarray}
with $E_0$ defined in (\ref{Ezero}) and $\w{K}_0$, $\w{K}_1$ defined in (\ref{K0K1}).
Setting $I_0$ and $I_1$ by integrals
\begin{equation}\label{i0}
I_0^2(j,|\gamma|) := \int_{\R^n}
\dfrac{|\partial_t^j(\lambda_{+}e^{\lambda_{-}t}
-\lambda_{-}e^{\lambda_{+}t})|^2}{|\lambda_{+}-\lambda_{-}|^2}\,|\xi|^{2|\gamma|}\chi(\xi)^2
d\xi
\end{equation}
and
\begin{equation}\label{i1}
I_1^2(j,|\gamma|) : = \int_{\R^n}
\dfrac{|\partial_t^j(e^{\lambda_{+}t}-e^{\lambda_{-}t})|^2}{|\lambda_{+}
-\lambda_{-}|^2}\,|\xi|^{2|\gamma|}\chi(\xi)^2 d\xi
\end{equation}
it follows from (\ref{K0K1}) that
\begin{equation}
||\partial_t^j \partial_x^{\gamma} E_0(t)||^2\lesssim
||v_0||_{L^1}^2I_0^2(j,|\gamma|)+||v_1||_{L^1}^2I_1^2(j,|\gamma|).\label{jalpha}
\end{equation}

To obtain the main results of this section, we estimate each one of
the integrals $I_0 $ and $ I_1$. To do this, we study separately two
cases, one where the eigenvalues are real and another where they are
complex. In both cases we use the lemma below to obtain the
estimates that we need.

\begin{lema}\label{integral}
Let $k>-n$, $\beta>0$ and $a>0$. Then
$$
\int_{|\xi| \leq \varepsilon} e^{-a|\xi|^{\beta}  t}|\xi|^kd\xi
\lesssim (1+t)^{-\frac{n+k}{\beta}},\quad \forall \,t>0.
$$
\end{lema}

%%%%%%%%%%%%%%%%%%%%%%%%%%%%%%%%%%%%%%%%%%%%%%%%%
%%%%%%%%%%%%%%%%%%%%%%%%%%%%%%%%%%%%%%%%%%%%%%%%%
%%%%%%%%% BAIXA FREQUENCIA CASO REAL %%%%%%%%%%%%
%%%%%%%%%%%%%%%%%%%%%%%%%%%%%%%%%%%%%%%%%%%%%%%%%
%%%%%%%%%%%%%%%%%%%%%%%%%%%%%%%%%%%%%%%%%%%%%%%%%

\vspace{0,5cm}\subsection{Case $\alpha>2\theta$: real eigenvalues}

\hspace{1cm}We consider $\theta\in [0,\frac{\alpha}{2} )$ and
$\varepsilon>0$ defined by
$\varepsilon^{\alpha-2\theta}=\frac{1}{4}$. For $|\xi|<\varepsilon$
we have $|\xi|^{2(\alpha-2\theta)}< \frac{1}{16}$ and thus the
eigenvalues are real. To estimate $I_0 $ and $I_1$ we use the
equivalences on the eigenvalues given by the following lemma.

\begin{lema}\label{autov} If $\,|\xi|<\varepsilon$ then:
\begin{enumerate}
\item[(i)] $\lambda_{+} \approx -|\xi|^{2(\alpha-\theta)}$ (more precisely,
$-4(2-\sqrt{2})|\xi|^{2(\alpha-\theta)}\leq\lambda_{+} \leq
-|\xi|^{2(\alpha-\theta)}$);
\item[(ii)] $\lambda_{-} \approx -|\xi|^{2\theta}$ (more precisely,
$-|\xi|^{2\theta}\leq\lambda_{-} \leq
-\frac{1}{4}\left(1+\frac{1}{\sqrt{2}}\right)|\xi|^{2\theta}$);
\item[(iii)] $\lambda_{+}-\lambda_{-}\approx  |\xi|^{2\theta}$.
\end{enumerate}
\end{lema}

\dem

\noindent \textit{(i)} First, we note that
$$|\xi|^{2(\alpha-2\theta)}< \frac{1}{16}=\frac{ 4\left(2-\sqrt{2}\,\right)-2}{16\left(2-\sqrt{2}\,\right)^2}.$$

Multiplying the above estimate by $4|\xi|^{2(\alpha-2\theta)}$ it
implies that
$$\begin{array}{ll}
-16\left(2-\sqrt{2}\,\right)|\xi|^{2(\alpha-2\theta)}+64\left(2-\sqrt{2}\,\right)^2|\xi|^{4(\alpha-2\theta)}
\leq -8|\xi|^{2(\alpha-2\theta)}\leq
-4|\xi|^{2(\alpha-2\theta)}\left(1+|\xi|^{2\delta}\right).
\end{array}$$

Using that
$1-16(2-\sqrt{2})|\xi|^{2(\alpha-2\theta)}+64(2-\sqrt{2})^2|\xi|^{4(\alpha-2\theta)}
=\left(1-8(2-\sqrt{2})|\xi|^{2(\alpha-2\theta)}\right)^2$ and
$1-8(2-\sqrt{2})|\xi|^{2(\alpha-2\theta)}>0$, we obtain
$$\begin{array}{ll}
-8\left(2-\sqrt{2}\,\right)|\xi|^{2(\alpha-2\theta)}\leq -1 +
\sqrt{1-4|\xi|^{2(\alpha-2\theta)}(1+|\xi|^{2\delta})}.
\end{array}$$

From $-\frac{1}{2}\leq -\frac{1}{2(1+|\xi|^{2\delta})}$ we conclude that
$$\begin{array}{ll}
-4\left(2-\sqrt{2}\,\right)|\xi|^{2(\alpha-\theta)}\leq
\displaystyle \frac{|\xi|^{2\theta}}{2(1+|\xi|^{2\delta})}
 \left(-1 + \sqrt{1-4|\xi|^{2(\alpha-2\theta)}\left(1+|\xi|^{2\delta}\right)}\right)= \lambda_{+}.
\end{array}$$

To prove that $\lambda_{+} \leq -|\xi|^{2(\alpha-\theta)}$ we see that
$$\left(-1 + \sqrt{1-4|\xi|^{2(\alpha-2\theta)}(1+|\xi|^{2\delta})}\right)
\leq -2|\xi|^{2(\alpha-2\theta)}\left(1+|\xi|^{2\delta}\right)$$
because $\,1-4
|\xi|^{2(\alpha-2\theta)}\left(1+|\xi|^{2\delta}\right)\leq
\left(1-2
|\xi|^{2(\alpha-2\theta)}\left(1+|\xi|^{2\delta}\right)\right)^2\,$
and $\,1-2|\xi|^{2(\alpha-2\theta)}(1+|\xi|^{2\delta})\geq 0$.
Multiplying both sides of above estimate by
$\frac{|\xi|^{2\theta}}{2(1+|\xi|^{2\delta})}$ we complete the
proof.

\vspace{0,3cm}\noindent \textit{(ii)} From $\frac{1}{2}\leq
1-8|\xi|^{2(\alpha-2\theta)}\leq
1-4|\xi|^{2(\alpha-2\theta)}(1+|\xi|^{2\delta})$ we can claim that
$$\frac{1}{4} \left(1+\frac{1}{\sqrt{2}}\right)\leq  \frac{1}{2(1+|\xi|^{2\delta})}
\left(1+\sqrt{1-4|\xi|^{2(\alpha-2\theta)}\left(1+|\xi|^{2\delta}\right)}\right)
\leq 1.$$ To conclude the proof we multiply both sides by
$-|\xi|^{2\theta}$.

\vspace{0,3cm}\noindent \textit{(iii)} By the choice of
$\varepsilon$, we have
$|\xi|^{2(\alpha-2\theta)}(1+|\xi|^{2\delta})\leq \frac{1}{8}$. Thus
$$\begin{array}{ll}
\displaystyle -\frac{1}{2}\leq
-4|\xi|^{2(\alpha-2\theta)}\left(1+|\xi|^{2\delta}\right) \, \leq 0.
\end{array}$$
Therefore
$$\frac{|\xi|^{2\theta}}{2\sqrt{2}}\leq \,
\frac{|\xi|^{2\theta}}{1+|\xi|^{2\delta}}
\,\sqrt{1-4|\xi|^{2(\alpha-2\theta)}\left(1+|\xi|^{2\delta}\right)}
\, \leq |\xi|^{2\theta}.$$

\fimdem

\begin{lema}\label{integrais} Let $I_0$ and $I_1$ given by {\rm(\ref{i0})} and {\rm(\ref{i1})},
$n\geq 1 $ and $\gamma$ multi-index. Then for all $t\geq 0$
we have:
\begin{enumerate}
\item[(i)] If $\,n+2|\gamma|> 4\theta$, then $I_1(0,|\gamma|)\lesssim
(1+t)^{-\frac{1}{\alpha-\theta}\left(\frac{n}{4}+\frac{|\gamma|}{2}-\theta\right)}$;
\item[(ii)] $I_1(1,|\gamma|)\lesssim
\begin{cases}
(1+t)^{-\frac{1}{\theta}\left(\frac{n}{4}
+\frac{|\gamma|}{2}\right)}+(1+t)^{-\frac{1}{\alpha-\theta}\left(\frac{n}{4}
+\frac{|\gamma|}{2}-\theta\right)-1}, & \text{ if } \,\theta>0 \\
(1+t)^{-\frac{1}{\alpha-\theta}\left(\frac{n}{4}
+\frac{|\gamma|}{2}-\theta\right)-1}, & \text{ if } \,\theta=0;
\end{cases}$
\item[(iii)] $I_0(0,|\gamma|)\lesssim (1+t)^{-\frac{1}{\alpha-\theta}\left(\frac{n}{4}
+\frac{|\gamma|}{2}\right)}$;
\item[(iv)] $I_0(1,|\gamma|)\lesssim (1+t)^{-\frac{1}{\alpha-\theta}\left(\frac{n}{4}
+\frac{|\gamma|}{2}\right)-1}$.
\end{enumerate}
\end{lema}

\dem

We will estimate (\ref{i1}) to obtain \textit{(i)} and \textit{(ii)}. The
items \textit{(iii)} and \textit{(iv)} can be proved in the same way using
(\ref{i0}).

\noindent \textit{(i)} From Lemma \ref{autov} follow that
$$\dfrac{|e^{\lambda_{+}t}-e^{\lambda_{-}t}|^2}{|\lambda_{+}-\lambda_{-}|^2}\lesssim
\dfrac{e^{-2|\xi|^{2(\alpha-\theta)}t}}{|\xi|^{4\theta}}.$$

Using the above estimate and the Lemma \ref{integral} we obtain
\begin{eqnarray*}
I_1^2(0,|\gamma|)&\lesssim & \int_{|\xi|<\varepsilon}
e^{-2|\xi|^{2(\alpha-\theta)}t}
|\xi|^{2|\gamma|-4\theta} d\xi  \\
&\lesssim & (1+t)^{-\frac{n+2|\gamma|-4\theta}{2(\alpha-\theta)}}
\end{eqnarray*}
because $n+2|\gamma|-4\theta>0$.

\newpage\vspace{0,2cm}\noindent \textit{(ii)} Newly, using the Lemma
\ref{autov},
$$\dfrac{|\partial_t(e^{\lambda_{+}t}-e^{\lambda_{-}t})|^2}{|\lambda_{+}-\lambda_{-}|^2}
\lesssim \dfrac{\lambda_+^2
e^{2\lambda_+t}+\lambda_-^2e^{2\lambda_-t}}{|\lambda_{+}-\lambda_{-}|^2}
\lesssim
|\xi|^{4(\alpha-2\theta)}e^{-2|\xi|^{2(\alpha-\theta)}t}+e^{-c|\xi|^{2\theta}t},$$
with $c=\frac{1}{2}\left(1+\frac{1}{\sqrt{2}}\right)$.

In the same way it was done in \textit{(i)}, we get for $\theta \in
(0,\frac{\alpha}{2})$
\begin{eqnarray*}
I_1^2(1,|\gamma|)&\lesssim &
\int_{|\xi|<\varepsilon}\left(|\xi|^{4(\alpha-2\theta)}
e^{-2|\xi|^{2(\alpha-\theta)}t}+e^{-c|\xi|^{2\theta}t}\right)|\xi|^{2|\gamma|} d\xi\\
& \lesssim
&(1+t)^{-\frac{n+2|\gamma|+4(\alpha-2\theta)}{2(\alpha-\theta)}}+
(1+t)^{-\frac{n+2|\gamma|}{2\theta}}.
\end{eqnarray*}

The case $\theta=0$ is immediate.

\fimdem

Using the previous lemma we prove the following result to the real
case for low frequency:

\begin{prop}\label{E_0real} Let $\alpha,\theta$ such that $0\leq 2\theta < \alpha$,
$\gamma$ a multi-index and $(v_0,\, v_1)\in \,
L^1(\R^n)\times L^1(\R^n)$. Then the following estimates are true,
for all $t\geq 0$:

\vspace{0,15cm}\noindent (i) If $\,n+2|\gamma|> 4\theta$ then
\begin{equation}\label{01>}
||\partial_x^\gamma  E_0(t)||\lesssim
||v_0||_{L^1}(1+t)^{-\frac{1}{\alpha-\theta}\left(\frac{n}{4}+\frac{|\gamma|}{2}\right)}
+||v_1||_{L^1}(1+t)^{-\frac{1}{\alpha-\theta}\left(\frac{n}{4}+\frac{|\gamma|}{2}-\theta\right)};
\end{equation}

\noindent (ii) If $\,n+2|\gamma|<2\alpha\,$ and $\,\theta\in
\left[\frac{n}{4}+\frac{|\gamma|}{2},\frac{\alpha}{2}\right)$ then
\begin{equation}\label{02ii>} ||\partial_t\partial_x^\gamma
E_0(t)||\lesssim
||v_0||_{L^1}(1+t)^{-\frac{1}{\alpha-\theta}\left(\frac{n}{4}+\frac{|\gamma|}{2}\right)-1}
+||v_1||_{L^1}(1+t)^{-\frac{1}{\theta}\left(\frac{n}{4}+\frac{|\gamma|}{2}\right)};
\end{equation}
and if $n+2|\gamma| \geq 2\alpha\,$ or $\,\theta\in
\left[0,\frac{n}{4}+\frac{|\gamma|}{2}\right) $ then
\begin{equation}\label{02i>}
||\partial_t\partial_x^\gamma E_0(t)||\lesssim
||v_0||_{L^1}(1+t)^{-\frac{1}{\alpha-\theta}\left(\frac{n}{4}+\frac{|\gamma|}{2}\right)-1}
+||v_1||_{L^1}(1+t)^{-\frac{1}{\alpha-\theta}\left(\frac{n}{4}+\frac{|\gamma|}{2}-\theta\right)-1}.
\end{equation}
\end{prop}

\dem

The item \textit{(i)} follows directly from items \textit{(i)} and
\textit{(iii)} of Lemma \ref{integrais} and from estimate
(\ref{jalpha}) with $j=0$. To prove the item \textit{(ii)}, we
choose $j=1$ in (\ref{jalpha}) and we use the Lemma \ref{integrais},
items \textit{(ii)} and \textit{(iv)}. To obtain the best decay rate
given in \textit{(ii)}, we see that $n+2|\gamma|<2\alpha$ and
$\theta\in
\left[\frac{n}{4}+\frac{|\gamma|}{2},\frac{\alpha}{2}\right)$
implies
$$\frac{n+2|\gamma|+4(\alpha-2\theta)}{2(\alpha-\theta)}\geq
\frac{n+2|\gamma|}{2\theta}.$$ Thus (\ref{02ii>}) is true. On the
other hand, if $n+2|\gamma| \geq 2\alpha$ or $\theta\in
\left[0,\frac{n}{4}+\frac{|\gamma|}{2}\right)$ we have the opposite
inequality. This implies the estimate (\ref{02i>}).

\fimdem

%%%%%%%%%%%%%%%%%%%%%%%%%%%%%%%%%%%%%%%%%%%%%%%%%
%%%%%%%%%%%%%%%%%%%%%%%%%%%%%%%%%%%%%%%%%%%%%%%%%
%%%%%%%%% BAIXA FREQUENCIA CASO COMPLEXO %%%%%%%%
%%%%%%%%%%%%%%%%%%%%%%%%%%%%%%%%%%%%%%%%%%%%%%%%%
%%%%%%%%%%%%%%%%%%%%%%%%%%%%%%%%%%%%%%%%%%%%%%%%%

\subsection{Case $\alpha\leq2\theta$: complex eigenvalues}

\hspace{1cm}Now we consider $\theta\in
\left[\frac{\alpha}{2},\alpha\right]\,$ and
$\,\varepsilon\in(0,1)$. In this case we have complex eigenvalues, given by
$$\lambda_{\pm}=\frac{|\xi|^{2\theta}}{2(1+|\xi|^{2\delta})}\left(-1
\pm
i\sqrt{4|\xi|^{2(\alpha-2\theta)}(1+|\xi|^{2\delta})-1}\right).$$
We obtain the results for complex case following the same steps of the real case.

\begin{lema}\label{autov_} If $\,|\xi|<\varepsilon$
then:
\begin{enumerate}
\item[(i)] $|\lambda_{+}-\lambda_{-}| \approx \, |\xi|^\alpha $;
\item[(ii)] $|\lambda_\pm|^2 \lesssim |\xi|^{2\alpha}$;
\item[(iii)] $|e^{\lambda_\pm t}|\lesssim e^{-\frac{1}{4}|\xi|^{2\theta}t}$.
\end{enumerate}
\end{lema}

\dem

\noindent \textit{(i)} Using $\,1 \leq 1+|\xi|^{2\delta} \leq 2\,$,
the proof follows directly from the inequality
$$\frac{1}{8}\left(4|\xi|^{2(\alpha-2\theta)}(1+|\xi|^{2\delta})-1\right)\,\leq\,
|\xi|^{2(\alpha-2\theta)} \,\leq\,
4|\xi|^{2(\alpha-2\theta)}(1+|\xi|^{2\delta})-1, \hspace{0,5cm} \forall \,\,0 < |\xi| < \varepsilon.$$

\noindent \textit{(iii)} By $\text{Re}(\lambda_\pm)\leq
-\frac{1}{4}|\xi|^{2\theta}$ we have $|e^{\lambda_\pm t}|\leq
e^{\text{Re}(\lambda_\pm)t}\leq e^{-\frac{1}{4}|\xi|^{2\theta}t}$.

\fimdem

The proof of the next lemma is quite similar to that of Lemma
\ref{integrais}. We omit the details.

\begin{lema}\label{integrais_} Let $I_0$ and $I_1$ given by {\rm(\ref{i0})} and {\rm(\ref{i1})},
$n\geq 1$ and $\gamma$ multi-index. Then, for all $t\geq 0$
we have for $\theta >0$:
\begin{enumerate}
\item[(i)] If $\,n+2|\gamma|> 2\alpha$, then $I_1(0,|\gamma|)\lesssim
(1+t)^{-\frac{1}{\theta}\left(\frac{n}{4}+\frac{|\gamma|}{2}-\frac{\alpha}{2}\right)}$;
\item[(ii)] $I_1(1,|\gamma|)\lesssim (1+t)^{-\frac{1}{\theta}\left(\frac{n}{4}
+\frac{|\gamma|}{2}\right)}$;
\item[(iii)] $I_0(0,|\gamma|)\lesssim (1+t)^{-\frac{1}{\theta}\left(\frac{n}{4}
+\frac{|\gamma|}{2}\right)}$;
\item[(iv)] $I_0(1,|\gamma|)\lesssim (1+t)^{-\frac{1}{\theta}\left(\frac{n}{4}
+\frac{|\gamma|}{2}+\frac{\alpha}{2}\right)}$.
\end{enumerate}
If $\alpha=\theta=0$ then it is immediate that $\,I_i(j,|\gamma|)
\lesssim e^{-\frac{t}{4}}$, $t\geq 0$, for all $i,j\in\{0,1\}$.
\end{lema}

From estimate (\ref{jalpha}) and Lemma \ref{integrais_}, we have the following result:

\newpage\begin{prop}\label{E_0comp} Let $\delta,\alpha \geq 0$, $\theta
\in \left[\frac{\alpha}{2},\alpha\right]$, $\gamma$ a multi-index
and $(v_0,\, v_1)\in L^1(\R^n)\times L^1(\R^n)$. Then, for $t\geq
0$, the following estimates are true:

\vspace{0,15cm}\noindent (i) If $\,n+2|\gamma|> 2\alpha$ and
$\theta>0$ then
$$||\partial_x^\gamma E_0(t)||\lesssim
||v_0||_{L^1}(1+t)^{-\frac{1}{\theta}\left(\frac{n}{4}+\frac{|\gamma|}{2}\right)}
+||v_1||_{L^1}(1+t)^{-\frac{1}{\theta}\left(\frac{n}{4}+\frac{|\gamma|}{2}
-\frac{\alpha}{2}\right)};$$

\noindent (ii) If $\,n\geq 1$ and $\theta>0$, then
$$||\partial_t\partial_x^\gamma E_0(t)||\lesssim
||v_0||_{L^1}(1+t)^{-\frac{1}{\theta}\left(\frac{n}{4}
+\frac{|\gamma|}{2}+\frac{\alpha}{2}\right)}+||v_1||_{L^1}(1+t)^{-
\frac{1}{\theta}\left(\frac{n}{4}+\frac{|\gamma|}{2}\right)};$$

\noindent (iii) If $\,\alpha=\theta=0$ then
$||\partial_t^j\partial_x^\gamma E_0(t)
||\lesssim\{||v_0||_{L^1}+||v_1||_{L^1}\}\,e^{-\frac{t}{4}}$, for
$j=0,1$.
\end{prop}

%%%%%%%%%%%%%%%%%%%%%%%%%%%%%%%%%%%%%%%%%%%%%%%%%
%%%%%%%%%%%%%%%%%%%%%%%%%%%%%%%%%%%%%%%%%%%%%%%%%
%%%%%%%%%%%%%%% ALTA FREQUENCIA %%%%%%%%%%%%%%%%%
%%%%%%%%%%%%%%%%%%%%%%%%%%%%%%%%%%%%%%%%%%%%%%%%%
%%%%%%%%%%%%%%%%%%%%%%%%%%%%%%%%%%%%%%%%%%%%%%%%%

\section{High frequency region: $|\xi|\geq \varepsilon$}

\hspace{1cm}The results obtained in the previous section complete
the estimates we need to prove the Theorems \ref{teo>} and
\ref{teo<} in the region of low frequency. Next, we are going to
obtain estimates for the region of high frequency and the required
regularity at initial data to obtain the desired decay. For this
purpose, we apply the multiplier method in Fourier space.

We consider, for $|\xi|\geq \varepsilon$ (with $\varepsilon\in(0,1)$
defined at previous section), the following auxiliary function:
\begin{equation}\label{rho}
\rho(\xi)=\begin{cases}\dfrac{\varepsilon^{2\alpha+2\delta-4\theta}
|\xi|^{2\theta}}{2(1+|\xi|^{2\delta})},& \text{ if
}\alpha+\delta\geq 2\theta
\vspace{3mm}\\
\dfrac{\varepsilon^{-2\alpha+4\theta}|\xi|^{2\alpha-2\theta}}{4},
&\text{ if }\alpha+\delta < 2\theta.
\end{cases}
\end{equation}
It is easy to prove that
\begin{equation}\label{rhomenor}
\rho(\xi)\leq \dfrac{|\xi|^{2\theta}}{2(1+|\xi|^{2\delta})} \quad
\text{and} \quad \rho(\xi) \leq \dfrac{|\xi|^{2\alpha-2\theta}}{2},
\end{equation}
for $|\xi| \geq \varepsilon$.

We denote by $E_1(t)$ the energy of order $\sigma$ of the equation
(\ref{f}) in the Fourier space given by
$$E_1(t)=\frac{1}{2}\,|\xi|^{2\delta+\sigma}
|\w{v}_t(t)|^2+\frac{1}{2}\,|\xi|^{2\alpha+
\sigma}|\w{v}(t)|^2,\quad t\geq 0.$$
Therefore
\begin{eqnarray}2\int_{|\xi|\geq \varepsilon }E_1(t)\,d\xi = \int_{\R^n}
(1-\chi(\xi))^2 \left\{|\xi|^{2\delta +\sigma}
|\w{v}_t(t)|^2+|\xi|^{2\alpha+\sigma} |\w{v}(t)|^2\right\}d\xi,\nonumber
\end{eqnarray}
where $\chi(\xi)$ is the characteristic function of $\{ \xi \in
\R^n\,\,/\,\,\, |\xi| < \varepsilon \}$.

Note that choosing $\sigma$ in a suitable way and using Plancherel
theorem, the above integral becames the $L^2$-norm of derivatives of
$E_\infty$ and $\partial_t E_\infty$ (see definition of $E_\infty$
at (\ref{Einf})). Thus our goal in this section is to estimate
$\int_{|\xi|\geq \varepsilon }E_1(t)d\xi $.

Next we show and prove the main result of this section. Note that
the decay estimate {\it (ii)} in the next proposition is due to the
regularity-loss structure of the equation for $\theta < \delta$. In
fact, to obtain decay rates to the energy in high frequency region
we need to impose more regularity on the initial data. The decay
rate is directly related to the additional regularity in the initial
data ($\frac{\delta-\theta}{\beta}$, with $\theta < \delta$). Thus,
to improve the decay rate it is necessary to take smaller $\beta$,
but this imply that it is necessary additional regularity on the
initial data. This does not occur in item {\it (i)}.

\begin{prop}\label{alta}
Let $\delta, \alpha \geq 0$, $\theta \in \left[0,\alpha\right]$
and $\sigma \in \R$.

\vspace{0,2cm}\noindent (i) If $(v_0,\, v_1)\in H^{\alpha+\frac{
\sigma}{2}}(\R^n)\times H^{\delta+\frac{ \sigma}{2}}(\R^n)$ and
$\delta\leq \theta \leq \alpha$, then there is a constant $c>0$ such
that
$$ \int_{|\xi|\geq \varepsilon} E_1(t,\xi)\,d\xi \,\lesssim
\, \left\{||v_1||_{H^{\delta+\frac{ \sigma}{2}}}^2
+ ||v_0||_{H^{\alpha+\frac{ \sigma}{2}}}^2\right\} e^{-ct},\quad
\forall \,t\geq 0.$$

\vspace{0,2cm}\noindent (ii) If $(v_0,\, v_1)\in H^{s}(\R^n)\times
H^{r}(\R^n)$, $\beta>0$ and $\theta < \delta$, then there is
$C=C(\beta)>0$ such that
$$\int_{|\xi|\geq \varepsilon} E_1(t,\xi)\,d\xi\,\leq C\left\{||v_1||_{H^r}^2+ ||v_0||_{H^{s}}^2\right\}
(1+t)^{-\frac{1}{\beta}},\quad \forall\, t\geq 0$$ with
$s=\alpha+\frac{\delta-\theta}{\beta} +\frac{ \sigma}{2}\, $ and
$\,r=\delta+\frac{\delta-\theta}{\beta} +\frac{ \sigma}{2}$.
\end{prop}

To prove Proposition \ref{alta} we need several estimates and lemmas
on the solution $\w{v}(t,\xi)$ of the corresponding problem {\rm
(\ref{f})} in the Fourier space. Multiplying equation (\ref{f}) by
$|\xi|^{\sigma}\over{\w{v}_t } + \rho(\xi)
|\xi|^{\sigma}\over{\w{v}}\,$ and taking the real part we have
\begin{equation}\label{EFR}
\frac{d}{dt}E(t)+F(t)=R(t),\quad t\geq 0
\end{equation}
where
\begin{eqnarray} \displaystyle
E(t,\xi)&=&\frac{1}{2}(1+|\xi|^{2\delta})|\xi|^{\sigma}|\w{v}_t(t)|^2
+\frac{1}{2}|\xi|^{2\alpha +\sigma}|\w{v}(t)|^2 \nonumber
\\
\displaystyle && +\rho(\xi)(1+|\xi|^{2\delta})
|\xi|^{\sigma}\text{Re}\{\w{v}_t(t)\over{\w{v}}(t)\}
+\frac{1}{2}\rho(\xi)|\xi|^{2\theta +\sigma}|\w{v}(t)|^2; \nonumber
\\
\displaystyle
F(t,\xi)&=&|\xi|^{2\theta+\sigma}|\w{v}_t(t)|^2+\rho(\xi)
|\xi|^{2\alpha +\sigma}|\w{v}(t)|^2; \nonumber
\\
\displaystyle
R(t,\xi)&=&\rho(\xi)(1+|\xi|^{2\delta})|\xi|^{\sigma}|\w{v}_t(t)|^2.\nonumber
\end{eqnarray}

By (\ref{rhomenor}) it follows that $R(t)\leq \frac{1}{2}F(t)$.
Substituting this estimate into (\ref{EFR}), we get
\begin{equation}\label{E-F}
\frac{d}{dt} E(t)+\frac{1}{2}F(t)\leq 0,
\end{equation}
for all $t\geq 0$ and $|\xi|\geq\varepsilon$.

\newpage\begin{lema}\label{equi} The functionals $E(t,\xi)$ and
$E_1(t,\xi)$ are equivalent for all $t\geq 0$ and
$|\xi|\geq\varepsilon$.
\end{lema}

\dem

First, we note that
\begin{eqnarray}\label{pm}
\pm
\rho(\xi)(1+|\xi|^{2\delta})\text{Re}\{\w{v}_t(t)\over{\w{v}}(t)\}&\leq &
(1+|\xi|^{2\delta})\frac{|\w{v}_t(t)|^2}{2\eta}+\rho(\xi)^2(1+|\xi|^{2\delta})
\frac{\eta|\w{v}(t)|^2}{2} \nonumber\\ &\leq &
(1+|\xi|^{2\delta})\frac{|\w{v}_t(t)|^2}{2\eta}+|\xi|^{2\alpha}
\frac{\eta|\w{v}(t)|^2}{8},\end{eqnarray} because, by (\ref{rhomenor}),
$$\rho(\xi)^2(1+|\xi|^{2\delta}) \leq\frac{|\xi|^{2\theta}}{2(1+|\xi|^{2\delta})}
\,\frac{|\xi|^{2\alpha-2\theta}}{2}\,(1+|\xi|^{2\delta})=
\dfrac{|\xi|^{2\alpha}}{4}.$$

Using (\ref{pm}) with positive sign and $\eta = 1$, we have:
\begin{eqnarray*}
E (t) \lesssim |\xi|^{\sigma} \left(
\frac{1}{2}|\xi|^{2\delta}|\w{v}_t(t)|^2
+\frac{1}{2}|\xi|^{2\alpha}|\w{v}(t)|^2\right)= E_1 (t) ,
\end{eqnarray*}
since, by (\ref{rhomenor}), $\rho(\xi)|\xi|^{2\theta}\leq
\frac{1}{2}|\xi|^{2\alpha}$ and $1+|\xi|^{2\delta}\leq
(\varepsilon^{-2\delta}+1)|\xi|^{2\delta}$.

On the other hand, using (\ref{pm}) with negative sign and $\eta =
2$, we have:
\begin{eqnarray*}
E (t) \geq |\xi|^{\sigma} \left( \frac{1}{4}|\xi|^{2\delta}
|\w{v}_t(t)|^2+\frac{1}{4}|\xi|^{2\alpha}|\w{v}(t)|^2
\right)=\frac{1}{2} \,E_1 (t).
\end{eqnarray*}

\fimdem

\begin{lema}\label{e1f}
If $\,\delta\leq \theta \leq \alpha$ then $E_1(t,\xi)\lesssim
F(t,\xi)$ for all $t\geq 0$ and $|\xi|\geq\varepsilon$.
\end{lema}

\dem

For all $t\geq 0\,$ and $\,|\xi|\geq\varepsilon$, we have
\begin{equation} \label{e1f1}|\xi|^{2\delta}|\w{v}_t(t)|^2
\lesssim |\xi|^{2\theta}|\w{v}_t(t)|^2\end{equation} and
\begin{equation} \label{e1f2} |\xi|^{2\alpha}|\w{v}(t)|^2
\lesssim \rho(\xi)|\xi|^{2\alpha}|\w{v}(t)|^2.\end{equation} The
last inequality actually works, because the case $\alpha+\delta \geq
2\theta\,$ follows from (\ref{rho}) using the hypothesis $\delta\leq
\theta$; and the case $\alpha+\delta <2\theta$ follows from
$|\xi|^{2\alpha}\gtrsim |\xi|^{2\theta}$, which works by hypothesis
$\theta\leq \alpha$.

Multiplying (\ref{e1f1}) and (\ref{e1f2}) by $|\xi|^{\sigma}$ and
adding the two  resulting inequalities, the lemma is proved.

\fimdem

Using the estimate (\ref{E-F}) and the Lemmas \ref{equi} and
\ref{e1f}, we can prove that there are positive constants $c, c_1$
such that
$$\frac{d}{dt} E(t)+c\,E(t)\,\leq \,\frac{d}{dt} E(t)+c_1 E_1(t)\,
\leq \, \frac{d}{dt} E(t)+\frac{1}{2} \, F(t) \leq 0, \hspace{0,5cm} \forall \,t \geq 0\,\,\,\text{and}\,\,\,|\xi| \geq \varepsilon.$$

By the equivalence between $E$ and $E_1$ given by Lemma \ref{equi}
we conclude that
$$E_1(t,\xi)\lesssim \, e^{-ct}E_1(0,\xi), \hspace{0,5cm} \forall \,t \geq 0\,\,\,\text{and}\,\,\,|\xi| \geq \varepsilon.$$

Integrating over the region of high frequency
$|\xi|\geq\varepsilon$, we obtain
\begin{eqnarray*}
\int_{|\xi|\geq \varepsilon} E_1(t,\xi)\,d\xi &\lesssim & e^{-ct} \int_{|\xi|\geq\varepsilon} \left\{
|\xi|^{2\delta+\sigma}
|\w{v}_1|^2 + |\xi|^{2\alpha+\sigma} |\w{v}_0|^2\right\}d\xi\\
& \lesssim & e^{-ct}
\int_{|\xi|\geq\varepsilon}(1+|\xi|^{2})^{\delta+\frac{\sigma}{2}}
|\w{v}_1|^2 d\xi +\int_{|\xi|\geq\varepsilon}
(1+|\xi|^2)^{\alpha+\frac{\sigma}{2}}
|\w{v}_0|^2d\xi\\
& \lesssim &  \left\{|| v_1||_{H^{\delta+\frac{\sigma}{2}}}^2 +
||v_0||_{H^{\alpha+\frac{\sigma}{2}}}^2\right\}e^{-ct}.
\end{eqnarray*}
The proof of the item \textit{(i)} of Proposition \ref{alta} is
complete.

The next lemma is important to estimate the integral of energy in
high frequency region in the Fourier space when there is the
property of regularity-loss, that is, when $\theta < \delta$. Using
the result of next lemma we will prove the item \textit{(ii)} of
Proposition \ref{alta}.

\begin{lema}\label{EF}
We define $$I(t)=\int_{|\xi|\geq \varepsilon} E_1(t,\xi)\,d\xi \quad
\text{and} \quad J(t)=\int_{|\xi|\geq \varepsilon} F(t,\xi)\,d\xi.$$
Let $\delta, \alpha \geq 0$, $\theta \in
\left[0,\alpha\right]$, $\sigma \in \R$, $\beta> 0$, $\theta <
\delta$ and $(v_0,\, v_1)\in H^{s}(\R^n)\times H^{r}(\R^n)$.
Then there is $C_\beta>0$ such that, for all $t\geq 0$,
$$ [I(t)]^{1+\beta}\leq \,C_\beta \{||v_1||_{H^r}+ ||v_0||_{H^{s}}\}^\beta J(t)$$
with $\,s=\alpha+\frac{\delta-\theta}{\beta} +\frac{ \sigma}{2}\,$ and
$\,r=\delta+\frac{\delta-\theta}{\beta} +\frac{ \sigma}{2}$.
\end{lema}

\dem

For any $\beta>0$ we have
\begin{eqnarray}
[I(t)]^{1+\beta} &\lesssim &\displaystyle \biggl[
\,\displaystyle{ \int_{|\xi|\geq \varepsilon}}
|\xi|^{2\delta+\sigma}\,|\w{v}_t(t)|^2\,d\xi\,\biggl]^{1+\beta} +
\,\biggl[ \displaystyle{ \int_{|\xi|\geq \varepsilon}}
|\xi|^{2\alpha+\sigma}\,|\w{v}(t)|^2\,d\xi\,\biggl]^{1+\beta}
\nonumber \\
&=& \displaystyle \biggl[ \,\displaystyle{
\int_{|\xi|\geq \varepsilon}}
|\xi|^{-\frac{2\theta+\sigma}{1+\beta}}
|\xi|^{2\delta+\sigma}\,|\w{v}_t(t)|^{2-\frac{2}{1+\beta}}\,
|\xi|^{\frac{2\theta+\sigma}{1+\beta}}
|\w{v}_t(t)|^{\frac{2}{1+\beta}}\,d\xi\,\biggl]^{1+\beta}
\nonumber \\
&&+ \displaystyle \biggl[ \displaystyle{ \int_{|\xi|\geq
\varepsilon}}
\Bigl(\rho(\xi)|\xi|^{2\alpha+\sigma}\Bigl)^{-\frac{1}{1+\beta}}|\xi|^{2\alpha+\sigma}
|\w{v}(t)|^{2-\frac{2}{1+\beta}}\,\Bigl(\rho(\xi)
|\xi|^{2\alpha+\sigma}\Bigl)^{\frac{1}{1+\beta}}
|\w{v}(t)|^{\frac{2}{1+\beta}}\,d\xi \biggl]^{1+\beta} \nonumber
\end{eqnarray}
for all $\,t \geq 0$. Then, by using H\"older's inequality in
$L^\frac{1+\beta}{\beta}$ and $L^{1+\beta}$ we obtain
\begin{equation}\label{As}
[I(t)]^{1+\beta} \lesssim \displaystyle \biggl\{ \biggl[ \,
\displaystyle{ \int_{|\xi|\geq \varepsilon}}
|\xi|^{-\frac{2\theta}{\beta}+\sigma+\frac{2\delta}{\beta}+2\delta}
|\w{v}_t(t)|^2\,d\xi\,\biggl]^{\beta} + \biggl[ \,\displaystyle{
\int_{|\xi|\geq \varepsilon}}
\rho(\xi)^{-\frac{1}{\beta}}|\xi|^{2\alpha+\sigma}\,
|\w{v}(t)|^2\,d\xi\,\biggl]^{\beta}\,\biggl\}\,J(t)
\end{equation}
for all $\,t \geq 0$.

Next we need to estimate both integrals in the right hand side of
the last inequality in terms of the initial data. Multiplying
equation (\ref{f}) by $\overline{\widehat{v}_t}$, taking the real
part on the resulting identity and integrating over the time
interval $(0,t)$ we obtain
$$\begin{array}{ll}
\dis (1+|\xi|^{2\delta})|\w{v}_t(t)|^2+|\xi|^{2\alpha}|\w{v}(t)|^2 +
2\int_0^t|\xi|^{2\theta}|\w{v}_t(s)|^2ds
=(1+|\xi|^{2\delta})|\w{v}_1|^2 +|\xi|^{2\alpha}|\w{v}_0|^2.
\end{array}$$

Now, multiplying both sides of the above identity by
$|\xi|^{-\frac{2\theta}{\beta}+\sigma+\frac{2\delta}{\beta}}$ and
integrating over $|\xi|\geq \varepsilon$ we obtain
\begin{eqnarray}
&&\displaystyle{\int_{|\xi|\geq
\varepsilon}|\xi|^{-\frac{2\theta}{\beta}
+\sigma+\frac{2\delta}{\beta}+2\delta}|\w{v}_t(t)|^2d\xi
+\int_{|\xi|\geq \varepsilon}|\xi|^{-\frac{2\theta}{\beta}+\sigma
+\frac{2\delta}{\beta}+2\alpha}|\w{v}(t)|^2 d\xi} \nonumber\\
&&\quad\displaystyle{\lesssim  \int_{|\xi|\geq
\varepsilon}(1+|\xi|^2)^{\delta+\frac{\delta-\theta}{\beta}
+\frac{\sigma}{2}}|\w{v}_1|^2d\xi+\int_{|\xi|\geq
\varepsilon}(1+|\xi|^{2})^{\alpha+\frac{\delta-\theta}{\beta}
+\frac{\sigma}{2}}|\w{v}_0|^2d\xi }. \label{sob}
\end{eqnarray}

By the assumptions of lemma we have $\alpha+\delta\geq 2
\theta$. Then, by definition of $\rho(\xi)$ given in (\ref{rho})
it follows that
\begin{equation}\label{A2>}
\int_{|\xi|\geq \varepsilon}
\rho(\xi)^{-\frac{1}{\beta}}|\xi|^{2\alpha+\sigma}\,
|\w{v}(t)|^2\,d\xi \leq K_\beta \int_{|\xi|\geq
\varepsilon}|\xi|^{-\frac{2\theta}{\beta} +\sigma
+\frac{2\delta}{\beta}+2\alpha}|\w{v}(t)|^2 d\xi.
\end{equation}

Using (\ref{sob})-(\ref{A2>}) in (\ref{As}), we obtain
$$[I(t)]^{1+\beta} \leq  C_\beta \left\{||v_1||_{H^r}^2+
||v_0||_{H^{s}}^2\right\}^\beta J(t), \quad \forall \,t\geq 0,$$
with $\,r=\delta+\frac{\delta-\theta}{\beta}+\frac{\sigma}{2}\,$ and
$\,s=\alpha+\frac{\delta-\theta}{\beta}+\frac{\sigma}{2}$.

\fimdem

Now, we observe that, by Lemma \ref{EF}, we can write
$$\dfrac{K}{2}\,[I(t)]^{1+\beta}\leq \dfrac{1}{2}\, J(t),\quad
\forall\, t\geq 0,$$ where $K=\left(C_\beta \{||v_1||_{H^r}^2+
||v_0||_{H^{s}}^2\}^\beta \right)^{-1}$.

By Lemma \ref{equi}, there are constants $m$ and $M$ such that
$$m E_1(t,\xi)\leq E(t,\xi) \leq ME_1(t,\xi),\quad \forall\, t\geq 0 \,\,\,\text{and}\,\,\,|\xi| \geq\varepsilon.$$
Consequently, we have
\begin{equation}\label{mM}
\frac{d}{dt}\int_{|\xi|\geq \varepsilon}E(t)\,d\xi+
\frac{K}{2M^{1+\beta}} \left[\int_{|\xi|\geq \varepsilon}E(t) \,d\xi
\right]^{1+\beta} \leq \frac{d}{dt}\int_{|\xi|\geq
\varepsilon}E(t)\,d\xi+ \frac{1}{2} \,J(t) \leq 0
\end{equation}
where the last inequality is due to the estimate (\ref{E-F}).

Therefore, applying the Lyapunov's theorem (see \cite{lya}) in
(\ref{mM}), we conclude the following decay estimate
$$\int_{|\xi|\geq \varepsilon} E_1(t,\xi)\,d\xi \leq C \left\{||v_1||_{H^r}^2+ ||v_0||_{H^{s}}^2 \right\}
(1+t)^{-\frac{1}{\beta}},\quad \forall\, t\geq 0,$$ with
$0<C=C(\beta)$. The proof of the item \textit{(ii)} of Proposition
\ref{alta} is complete.

%%%%%%%%%%%%%%%%%%%%%%%%%%%%%%%%%%%%%%%%%%%%%%%%%
%%%%%%%%%%%%%%%%%%%%%%%%%%%%%%%%%%%%%%%%%%%%%%%%%
%%%%%%%%%%%% RESULTADOS PRINCIPAIS %%%%%%%%%%%%%%
%%%%%%%%%%%%%%%%%%%%%%%%%%%%%%%%%%%%%%%%%%%%%%%%%
%%%%%%%%%%%%%%%%%%%%%%%%%%%%%%%%%%%%%%%%%%%%%%%%%

\section{Results and applications}

\subsection{Main results}

\hspace{1cm}In this section we use the results obtained in the
previous sections in order to find decay rates for the $L^2$-norm of
$\partial_x^{\gamma_1}v(t)$ and $\partial_x^{\gamma_2}v_t(t)$ and
the corresponding regularity in the initial data. Observe that for
appropriate choices of $\gamma_1$ and $\gamma_2$ we can obtain
estimates for the $L^2$-norm of solution, for each term of energy
defined in (\ref{energia}) and also for $H^m$ norms of $v$ and
$v_t$. In Theorem \ref{teo>} we do not determine values for
$\beta>0$ for the purpose of to make this choice appropriately on
each application. The same occurs in Theorem \ref{teo<}. We can
choose $\beta$ so that give the optimal rates and in this case will
require more of the regularity of the initial data. Or we can also
choose $\beta$ to require less regularity on the initial data,
resulting in worse rates.

\begin{teo}\label{teo>} If $\delta,\,\alpha \geq 0$,
$\theta\in \left[0,\frac{\alpha}{2}\right)$, $v_0 \in
H^{s}(\R^n)\cap\, L^1(\R^n)$ and $v_1 \in H^r(\R^n)\cap\,
L^1(\R^n)$, for $r$ and $s$ specified below, then the solution of the
problem {\rm(\ref{eq})-(\ref{dados})} satisfies the following
estimates for all $t \geq 0$:

\vspace{0,2cm}\noindent (i) If $\,n+2|\gamma_1|>4\theta$ then
\begin{eqnarray}\label{estimav}
||\partial_x^{\gamma_1}v(t)||&\lesssim
&||v_0||_{L^1}(1+t)^{-\frac{1}{\alpha-\theta}\left(\frac{n}{4}+\frac{|\gamma_1|}{2}\right)}
+||v_1||_{L^1}(1+t)^{-\frac{1}{\alpha-\theta}\left(\frac{n}{4}
+\frac{|\gamma_1|}{2}-\theta\right)}\nonumber
\\ &&+\,\left\{||v_1||_{H^r}+
||v_0||_{H^{s}}\right\}\,f(t),\nonumber
\end{eqnarray} with $f(t) = \left\{
\begin{array}{ll} e^{-ct}\,\,(c \in \R) \hspace{0,37cm}
\text{if} \,\,\,\delta \leq \theta  \\
(1+t)^{-\frac{1}{2\beta}} \hspace{0,7cm} \text{if} \,\,\,\theta <
\delta \end{array} \right.$, $\,\,s = \left\{
\begin{array}{ll} |\gamma_1| \hspace{1,61cm}
\text{if} \,\,\,\delta \leq \theta  \\
|\gamma_1| + \frac{\delta-\theta}{\beta} \hspace{0,5cm} \text{if}
\,\,\,\theta < \delta \end{array} \right.$ and $\,r=s+ \delta
-\alpha$.

\vspace{0,3cm} \noindent (ii) If $\,n+2|\gamma_2|<2\alpha\,$ and
$\,\theta\in
\left[\frac{n}{4}+\frac{|\gamma_2|}{2},\frac{\alpha}{2}\right)$
then
\begin{eqnarray}\label{estivavt1}
||\partial_x^{\gamma_2}v_t(t)||&\lesssim &
||v_0||_{L^1}(1+t)^{-\frac{1}{\alpha-\theta}\left(\frac{n}{4}+\frac{|\gamma_2|}{2}\right)-1}
+||v_1||_{L^1}(1+t)^{-\frac{1}{\theta}\left(\frac{n}{4}+\frac{|\gamma_2|}{2}\right)}\nonumber\\
&&+\left\{||v_1||_{H^r}+ ||v_0||_{H^{s}}\right\}f(t),
\end{eqnarray}
and if $\,n+2|\gamma_2| \geq 2\alpha\,$ or $\,\theta\in
\left[0,\frac{n}{4}+\frac{|\gamma_2|}{2}\right)$ then
\begin{eqnarray}\label{estimavt2}
||\partial_x^{\gamma_2} v_t(t)||& \lesssim  &
||v_0||_{L^1}(1+t)^{-\frac{1}{\alpha-\theta}\left(\frac{n}{4}+\frac{|\gamma_2|}{2}\right)-1}
+||v_1||_{L^1}(1+t)^{-\frac{1}{\alpha-\theta}\left(\frac{n}{4}
+\frac{|\gamma_2|}{2}-\theta\right)-1}\nonumber \\
&&+\left\{||v_1||_{H^r}+ ||v_0||_{H^{s}}\right\}f(t),
\end{eqnarray} with $f(t) = \left\{
\begin{array}{ll} e^{-ct}\,\,(c \in \R) \hspace{0,37cm}
\text{if} \,\,\,\delta \leq \theta  \\
(1+t)^{-\frac{1}{2\beta}} \hspace{0,7cm} \text{if} \,\,\,\theta <
\delta \end{array} \right.$, $\,\,r = \left\{
\begin{array}{ll} |\gamma_2 |\hspace{1,61cm}
\text{if} \,\,\,\delta \leq \theta \\
|\gamma_2| + \frac{\delta-\theta}{\beta} \hspace{0,5cm} \text{if}
\,\,\,\theta < \delta \end{array} \right.$ and $\,s=r + \alpha
-\delta$.
\end{teo}

\dem

\noindent \textit{(i)} By definition of $E_0$ and $E_\infty$ we have
$$
||\partial_x^{\gamma_1} v(t)|| \leq  ||\partial_x^{\gamma_1}
E_0(t)|| + ||\partial_x^{\gamma_1} E_\infty(t)||.
$$
We choose $\gamma = \gamma_1 $ at Proposition \ref{E_0real} and
consider $n+2|\gamma_1|>4\theta$. Using (\ref{01>})
and the Proposition \ref{alta} with $\sigma = 2|\gamma_1| - 2
\alpha\,$ it follows the result.

\vspace{0,2cm}\noindent \textit{(ii)} We can write
$$
||\partial_x^{\gamma_2} v_t(t)|| \leq
||\partial_x^{\gamma_2}\partial_t E_0(t)|| +
||\partial_t\partial_x^{\gamma_2} E_\infty(t)||.
$$
Choosing $\gamma = \gamma_2$ at Proposition \ref{E_0real} and
$\sigma = 2|\gamma_2| - 2 \delta$ at Proposition \ref{alta}, we have
two options. If $n+2|\gamma_2|<2\alpha\,$ and $\,\theta\in
\left[\frac{n}{4}+\frac{|\gamma_2|}{2},\frac{\alpha}{2}\right)$ we
use (\ref{02ii>}) and obtain (\ref{estivavt1}). For the case
$n+2|\gamma_2| \geq 2\alpha\,$ or $\,\theta\in
\left[0,\frac{n}{4}+\frac{|\gamma_2|}{2}\right)$ using (\ref{02i>})
it follows (\ref{estimavt2}).

\fimdem

In the same way, we prove the following result, for the case $ \alpha\leq 2\theta$:

\begin{teo}\label{teo<} If $\delta,\,\alpha \geq 0$,
$\theta\in \left[\frac{\alpha}{2}, \alpha \right]$, $v_0 \in
H^{s}(\R^n)\cap\, L^1(\R^n)$ and $v_1 \in H^r(\R^n)\cap\,
L^1(\R^n)$, with $f(t)$, $r$ and $s$ defined at same manner as in
Theorem {\rm\ref{teo>}}, the following estimates are true, for all
$t \geq 0$:

\vspace{0,2cm}\noindent (i) If $\,n+2|\gamma_1|> 2\alpha$ then for
$\theta>0$ we have
$$||\partial_x^{\gamma_1} v(t)||\lesssim
||v_0||_{L^1}(1+t)^{-\frac{1}{\theta}\left(\frac{n}{4}+\frac{|\gamma_1|}{2}\right)}
+||v_1||_{L^1}(1+t)^{-\frac{1}{\theta}\left(\frac{n}{4}+\frac{|\gamma_1|}{2}
-\frac{\alpha}{2}\right)}+\left\{||v_1||_{H^r}+
||v_0||_{H^{s}}\right\}f(t),$$ and, for $\theta=0$,
$$||\partial_x^{\gamma_1}v(t)|| \lesssim
\{||v_0||_{L^1}+||v_1||_{L^1}\}\,e^{-\frac{t}{4}}+\left\{||v_1||_{H^r}+
||v_0||_{H^{s}}\right\}f(t);$$

\noindent (ii) If $\,n\geq 1$ and $\theta>0$ then
$$||\partial_x^{\gamma_2} v_t(t)||\lesssim
||v_0||_{L^1}(1+t)^{-\frac{1}{\theta}\left(\frac{n}{4}
+\frac{|\gamma_2|}{2}+\frac{\alpha}{2}\right)}+||v_1||_{L^1}(1+t)^{-
\frac{1}{\theta}\left(\frac{n}{4}+\frac{|\gamma_2|}{2}\right)}
+\left\{||v_1||_{H^r}+ ||v_0||_{H^{s}}\right\}f(t),$$ and, if
$\,\theta=0$, $$||\partial_x^{\gamma_2}v_t(t)||\lesssim
\{||v_0||_{L^1}
+||v_1||_{L^1}\}\,e^{-\frac{t}{4}}+\left\{||v_1||_{H^r}+
||v_0||_{H^{s}}\right\}f(t).$$
\end{teo}

%%%%%%%%%%%%%%%%%%%%%%%%%%%%%%%%%%%%%%%%%%%%%%%%%
%%%%%%%%%%%%%%%%%%%%%%%%%%%%%%%%%%%%%%%%%%%%%%%%%
%%%%%%%%%%%%%%%% APLICAÇÕES %%%%%%%%%%%%%%%%%%%%%
%%%%%%%%%%%%%%%%%%%%%%%%%%%%%%%%%%%%%%%%%%%%%%%%%
%%%%%%%%%%%%%%%%%%%%%%%%%%%%%%%%%%%%%%%%%%%%%%%%%

\subsection{Applications}

\hspace{1cm}Now let us apply the above results to several initial value problems associated with
some dissipative partial differential equations of second order in
time.

\subsubsection{Wave equation with fractional damping}

\hspace{1cm}We consider the equation (\ref{eq})-(\ref{dados}) with
$\delta =0\,$ and $\,\alpha =1 $, that is, the wave equation
$$v_{tt}(t,x) -\Delta v(t,x) + (-\Delta)^{\theta}v_t (t,x)= 0, \quad
t\geq 0,\,\,\, x\in \R^n$$
with initial data
\begin{equation*}
v(0,x) = v_0(x), \quad v_t(0,x)=v_1(x),
\end{equation*}
and $\theta \in [0,1]$, where the associated energy is defined by
\begin{equation*}E_v(t) = \frac{1}{2}\left\{ ||v_t(t)||^2 +  ||\nabla v(t)||^2 \right\}.\end{equation*}

Note that $\delta\leq \theta$. Moreover, for this choice of $\alpha$, we have $\alpha > 2\theta$ if and only
if $\theta<\frac{1}{2}$. Thus, it appears here the separation in
cases $\theta\in \left[0,\frac{1}{2}\right)$ (real eigenvalues for the low frequency region) and
$\theta \in \left[\frac{1}{2},1\right]$ (complex eigenvalues for the low frequency region). If
$\theta\in \left[0,\frac{1}{2}\right)$, we can apply Theorem
\ref{teo>}, to get that
$$
||v(t)||\lesssim \left\{||v_1||_{H^{-1}\cap L^1}+ ||v_0||_{L^2\cap
L^1}\right\}
(1+t)^{-\frac{1}{1-\theta}\left(\frac{n}{4}-\theta\right)},\quad \text{for all}\,\,n>4\theta;
$$ 
$$
||v_t(t)||\lesssim \begin{cases}\left\{||v_1||_{L^2\cap L^1} +
||v_0||_{H^1\cap L^1}\right\} (1+t)^{-\frac{1}{4\theta}},
& \text{if}\,\,\, n=1\,\,\, \text{and}\,\,\, \theta\in \left[\frac{1}{4},\frac{1}{2}\right),\\
\left\{||v_1||_{L^2\cap L^1}+ ||v_0||_{H^1\cap L^1}\right\}
(1+t)^{-\frac{1}{1-\theta}\left(\frac{n}{4}-\theta\right)-1},&
\text{if}\,\,\,  n\geq 2\,\,\, \text{or} \,\,\,\theta\in
\left[0,\frac{1}{4}\right),\end{cases}
$$
$$
||\nabla v(t)||\lesssim \left\{||v_1||_{L^2\cap L^1}+
||v_0||_{H^1\cap
L^1}\right\}(1+t)^{-\frac{1}{1-\theta}\left(\frac{n}{4}+\frac{1}{2}-\theta\right)},\quad \text{for}\,\, n\geq 1.
$$

Now for $\theta\in \left[\frac{1}{2},1\right]$ we have $\alpha\leq
2\theta$. Then, using Theorem \ref{teo<} we can obtain the following
estimate for the norm of solution, for $n\geq 3$:
$$
||v(t)||\lesssim \left\{||v_1||_{H^{-1}\cap L^1}+ ||v_0||_{L^2\cap
L^1}\right\}(1+t)^{-\frac{1}{\theta}\left(\frac{n}{4}-\frac{1}{2}\right)};
$$
and the following decay rates for the energy terms
$$
||v_t(t)||+||\nabla v(t)||\lesssim \left\{||v_1||_{L^2\cap L^1}+
||v_0||_{H^1\cap L^1} \right\} (1+t)^{-\frac{n}{4\theta}},
$$
in the case $n \geq 1$.

\subsubsection{Plate equation with rotational inertia effects and fractional damping}

\hspace{1cm}In the previous works \cite{RCR2, DRC1} for the plate
equation with rotational inertia effects and fractional damping the
authors considered the hypothesis $\theta \in [0,1]$. In this
application we generalize this interval assuming that $\theta \in
[0,2]$. Note that $\alpha>2\theta$ if $\theta\in[0,1)$ and $\alpha
\leq 2\theta$ if $\theta\in[1,2]$, that is, to obtain the results it
will be necessary to separate in two cases. This situation is
similar to that one for the wave equation with fractional damping
when the interval $\theta \in [0,1]$ is separated at intervals
$\theta\in\left[0,\frac{1}{2}\right)$ and $\theta\in
\left[\frac{1}{2},1\right]$ (see \cite{RCR1, IN, DR}).

We consider the following problem:
\begin{equation*}
\left\{\begin{array}{ll} v_{tt}(t,x)-\Delta v_{tt}(t,x)+\Delta^2
v(t,x) + (-\Delta)^{\theta}v_t (t,x)= 0, \quad t\geq 0,\,\,\, x\in
\R^n
\\
v(0,x) = v_0(x), \quad v_t(0,x)=v_1(x),
\end{array}\right.
\end{equation*}
with $\theta \in [0,2]$. That is, we take $\delta=1$ and $\alpha =
2$ in (\ref{eq})-(\ref{dados}).

The associated energy is given by
\begin{equation*}\label{energia-placas}E_v(t) = \frac{1}{2}\left\{ ||v_t(t)||^2
+ ||\nabla v_t(t)||^2 + ||\Delta v(t)||^2 \right\}.\end{equation*}

If $\theta\in[0,1)$, we consider $|\gamma_1|=0\,$ and
$\,\frac{1}{2\beta} =
\frac{1}{2-\theta}\left(\frac{n}{4}-\theta\right)$. By item
\textit{(i)} from Theorem \ref{teo>}, we have the following estimate
for the $L^2$-norm of solution, for all $n>4\theta$:
$$
||v(t)||\lesssim \left\{||v_1||_{H^r\cap L^1}+ ||v_0||_{H^{s}\cap
L^1}\right\}(1+t)^{-\frac{1}{2-\theta}\left(\frac{n}{4}-\theta\right)},
$$
with $\,r=\frac{(1-\theta)(n-4\theta)}{2(2-\theta)} -1\,$ and
$\,s=\frac{(1-\theta)(n-4\theta)}{2(2-\theta)}$.

In order to get decay rates for the energy, we estimate each term of
it separately. In the next three estimates we assume
$\theta\in[0,1)$. At item \textit{(ii)} of Theorem \ref{teo>} we can
choose $|\gamma_2|=0$ with $\frac{1}{2\beta}=\frac{n}{4\theta}$ in
(\ref{estivavt1}) and
$\frac{1}{2\beta}=\frac{1}{2-\theta}\left(\frac{n}{4}-\theta\right)+1$
in (\ref{estimavt2}), to obtain
$$%\begin{equation}\label{vt1-placas}
||v_t(t)||\lesssim \begin{cases}\left\{||v_1||_{H^{r_1}\cap L^1} +
||v_0||_{H^{s_1}\cap L^1}\right\}(1+t)^{-\frac{n}{4\theta}},
& \text{if}\,\,\, n<4\,\, \text{and}\,\, \theta\in \left[\frac{n}{4},1\right),\\
\left\{||v_1||_{H^{r_2}\cap L^1}+ ||v_0||_{H^{s_2}\cap L^1}\right\}
(1+t)^{-\frac{1}{2-\theta}\left(\frac{n}{4}-\theta\right)-1},&
\text{if}\,\,\,  n\geq 4\,\, \text{or} \,\,\theta\in
\left[0,\frac{n}{4}\right),\end{cases}
$$%\end{equation}
with ${r_1}=\frac{(1-\theta)n}{2\theta},\,$
${s_1}=\frac{(1-\theta)n}{2\theta}+1,\,$
$r_2=\frac{(1-\theta)(n-8\theta+8)}{2(2-\theta)}\,$ and
$\,s_2=\frac{(1-\theta)(n-8\theta+8)}{2(2-\theta)}+1$.

\vspace{0,1cm}Now, if $|\gamma_2|=1$, with
$\frac{1}{2\beta}=\frac{1}{\theta}\left(\frac{n}{4}+\frac{1}{2}\right)$
in (\ref{estivavt1}) and
$\frac{1}{2\beta}=\frac{1}{2-\theta}\left(\frac{n}{4}+\frac{1}{2}-\theta\right)+1$
in (\ref{estimavt2}), by item \textit{(ii)} in Theorem \ref{teo>} it
follows that
$$
||\nabla v_t(t)||\lesssim \begin{cases}\left\{||v_1||_{H^{r_1}\cap L^1}
+ ||v_0||_{H^{s_1}\cap
L^1}\right\}(1+t)^{-\frac{3}{4\theta}},& \text{if}\,\,\, n=1\,\,
\text{and}\,\, \theta\in \left[\frac{3}{4},1\right),\\
\left\{||v_1||_{H^{r_2}\cap L^1}+ ||v_0||_{H^{s_2}\cap
L^1}\right\}(1+t)^{-\frac{1}{2-\theta}\left(\frac{n}{4}+\frac{1}{2}-\theta\right)-1},&
\text{if}\,\,\, n\geq 2 \,\, \text{or} \,\, \theta\in
\left[0,\frac{3}{4}\right),
\end{cases}
$$
where $r_1=1+\frac{(1-\theta)(n+2)}{2\theta} $,
$s_1=2+\frac{(1-\theta)(n+2)}{2\theta}$,
$r_2=1+\frac{(1-\theta)(n-8\theta+10)}{2(2-\theta)}\,$ and
$\,s_2=2+\frac{(1-\theta)(n-8\theta+10)}{2(2-\theta)} $.

\vspace{0,1cm}On the other hand, if $|\gamma_1|=2$ and $\frac{1}{2\beta}
=\frac{1}{2-\theta}\left(\frac{n}{4}+1-\theta\right) $, the Theorem
\ref{teo>}, item \textit{(i)}, gives us:
$$
%\begin{equation}\label{lap-placas}
||\Delta v(t)||\lesssim \left\{||v_1||_{H^r\cap L^1}+
||v_0||_{H^s\cap
L^1}\right\}(1+t)^{-\frac{1}{2-\theta}\left(\frac{n}{4}+1-\theta\right)}
%\end{equation}
$$
where $\,r=1+\frac{(1-\theta)(n-4\theta+4)}{2(2-\theta)}\,$ and
$\,s=2+\frac{(1-\theta)(n-4\theta+4)}{2(2-\theta)}$.

\vspace{0,1cm}The result for the case $\theta \in[1,2]$ is obtained
similarly to the previous case, using the Theorem \ref{teo<} instead
of Theorem \ref{teo>}. The decay rates obtained to the $L^2$-norm of
solution are:
$$\begin{array}{ll}
||v(t)||\lesssim \left\{||v_1||_{H^{-1}\cap L^1}+ ||v_0||_{L^2\cap
L^1}\right\}(1+t)^{-\frac{1}{\theta}\left(\frac{n}{4}-1\right)},
\end{array}$$
for $n>4$. For the energy terms we have:
$$\begin{array}{ll}
||v_t(t)||\lesssim ||v_t(t)||\lesssim \left\{||v_1||_{L^2\cap L^1}+
||v_0||_{H^1\cap L^1} \right\}(1+t)^{-\frac{n}{4\theta}};
\\ \vspace{-0,3cm} \\
||\nabla v_t(t)||\lesssim \left\{||v_1||_{H^1\cap L^1}+
||v_0||_{H^2\cap L^1} \right\}
(1+t)^{-\frac{1}{\theta}\left(\frac{n}{4}+\frac{1}{2}\right)};
\\ \vspace{-0,3cm} \\
||\Delta v(t)||\lesssim \left\{||v_1||_{H^1\cap L^1}+
||v_0||_{H^2\cap L^1}\right\} (1+t)^{-\frac{n}{4\theta}}.
\end{array}$$

\begin{obs}
The plate equation without the inertia rotational effects, that is,
{\rm(\ref{eq})-(\ref{dados})} with $\alpha=2$ and $\delta = 0$, does
not have the structure of regularity-loss. In this case, the decay
rates are equal to the above rates, but without assume additional
regularity at the initial data.
\end{obs}

\subsubsection{Boussinesq equation with fractional damping}

\hspace{1,2cm}In this subsection, we want to show that it is
possible to add terms of type $(-\Delta)^{\delta_1}v_{tt}$,
$(-\Delta)^{\theta_1}v_t$, $(-\Delta)^{\alpha_1}v$ in the equation
(\ref{eq}) in order to obtain decay rates for the total energy and
for the $L^2$-norm of solution from the previous results. However,
we consider just a simple example to illustrate this case.

In \cite{PE} the authors studied the Boussinesq equation (IBq) with
the strong dissipation $\Delta v_t$ (see \cite{wang2,wang}). Next,
we regard a more general case than the linear problem studied in
\cite{PE}, considering the following IBq with fractional damping in
$\R^n$:
\begin{equation}\label{ibq}v_{tt}(t,x)-\Delta v(t,x)-\Delta v_{tt}(t,x)
+\Delta^2 v(t,x) + (-\Delta)^{\theta}v_t (t,x)= 0, \quad t\geq
0,\,\,\, x\in \R^n
\end{equation}
with $\theta\in\left[0,1\right]$ and initial data
\begin{equation}\label{ibqdados}
v(0,x) = v_0(x), \quad v_t(0,x)=v_1(x).
\end{equation}

The energy associated to this equation is given by
$$E_v(t)=\frac{1}{2}\left\{ ||v_t(t)||^2+||\nabla v (t)||^2
+||\nabla v_t(t)||^2+||\Delta v(t)||^2 \right\}.$$

Applying the Fourier transform in (\ref{ibq}) and
(\ref{ibqdados}), we obtain
$$\left\{
\begin{array}{ll}(1+|\xi|^{2})\widehat{v}_{tt}(t,\xi)+
|\xi|^{2\theta}\widehat{v}_t(t,\xi)+\left(|\xi|^{2}+|\xi|^{4}\right)\widehat{v}(t,\xi)=0,
\\ \widehat{v}(0,\xi)=\widehat{v}_0(\xi), \quad
\widehat{v_t}(0,\xi)=\widehat{v}_1(\xi).
\end{array}\right.$$

This problem is not a specific case of the initial problem
(\ref{f}). But we note that for $|\xi|<\varepsilon$ is valid the
equivalence
$$|\xi|^2 \approx |\xi|^2+|\xi|^4,$$
which leads us to conclude that the estimates in the low frequency
for this equation are the same of Propositions \ref{E_0real} and \ref{E_0comp} considering $\alpha=\delta=1$.

On the other hand, if $|\xi|\geq \varepsilon$, then
$$|\xi|^4 \approx |\xi|^2+|\xi|^4,$$
and therefore the estimates and regularity in the initial data for this equation
at high frequency are given by Proposition \ref{alta}, with $\alpha=2$ and $\delta=1$. This gives us the following results:

\begin{teo}\label{ibqteo} If $\,\theta\in \left[0,\frac{1}{2}\right)$ and
$(v_0,\, v_1)\in [H^{s}(\R^n)\cap\, L^1(\R^n)]\times
[H^r(\R^n)\cap\, L^1(\R^n)]$, with $r$ and $s$ specified in each
case below, then the solution of the problem
{\rm(\ref{ibq})-(\ref{ibqdados})} satisfies the following decay
estimates:

\vspace{0,15cm} \noindent (i)  If $\,n>4\theta$ then for $t\geq 0$
$$\begin{array}{ll}
||v(t)|| \lesssim \left\{||v_1||_{H^r\cap L^1}+||v_0||_{H^s\cap
L^1}\right\}(1+t)^{-\frac{1}{1-\theta}\left(\frac{n}{4}-\theta\right)},\nonumber
\end{array}$$
with $\,r=\frac{n-4\theta}{2}-1\,$ and $\,s = \frac{n-4\theta}{2}$.

\vspace{0,15cm} \noindent (ii)  If $\,n\geq 1$ then for $t\geq 0$
$$\begin{array}{ll}
||\nabla v(t)|| \lesssim \left\{||v_1||_{H^r\cap
L^1}+||v_0||_{H^s\cap L^1}\right\}
(1+t)^{-\frac{1}{1-\theta}\left(\frac{n}{4}+\frac{1}{2}-\theta\right)},\nonumber
\end{array}$$
with $\,r=\frac{n+2-4\theta}{2}\,$ and $\,s = 1 + \frac{n+2-4\theta}{2}
$;
$$\begin{array}{ll}
||\Delta v(t)|| \lesssim \left\{||v_1||_{H^r\cap
L^1}+||v_0||_{H^s\cap
L^1}\right\}(1+t)^{-\frac{1}{1-\theta}\left(\frac{n}{4}+1-\theta\right)},\nonumber
\end{array}$$
with $\,r=1 + \frac{n+4-4\theta}{2}\,$ and $\,s = 2 +
\frac{n+4-4\theta}{2}$;
$$\begin{array}{ll}
||v_t(t)|| \lesssim \begin{cases} \left\{||v_1||_{H^{r_1}\cap L^1}+
||v_0||_{H^{s_1}\cap L^1}\right\} (1+t)^{-\frac{1}{4\theta}}, &
\text{if}\,\,\, n=1\,\,\, \text{and}\,\,\, \theta\in
\left[\frac{1}{4},\frac{1}{2}\right),\\
\left\{||v_1||_{H^{r_2}\cap L^1}+ ||v_0||_{H^{s_2}\cap L^1}\right\}
(1+t)^{-\frac{1}{1-\theta}\left(\frac{n}{4}-\theta\right)-1},&
\text{if}\,\,\,  n\geq 2\,\,\, \text{or} \,\,\,\theta\in
\left[0,\frac{1}{4}\right),\end{cases}\nonumber
\end{array}$$
with $\,r_1=\frac{(1-\theta)n}{2\theta} $,
$s_1=1+\frac{(1-\theta)n}{2\theta}$, $r_2=\frac{n-8\theta+4}{2}\,$ and
$\,s_2=1+\frac{n-8\theta+4}{2} $;
$$\begin{array}{ll}
||\nabla v_t(t)|| \lesssim \left\{||v_1||_{H^r\cap L^1}+
||v_0||_{H^s\cap
L^1}\right\}(1+t)^{-\frac{1}{1-\theta}\left(\frac{n}{4}+\frac{1}{2}-\theta\right)-1},
\end{array}$$
with $\,r=1+\frac{n-8\theta+6}{2}\,$ and $\,s=2+\frac{n-8\theta+6}{2}
$.
\end{teo}

\begin{teo} If $\,\theta\in\left[\frac{1}{2}, 1\right]$ and $(v_0,\, v_1)\in
[H^{s}(\R^n)\cap\, L^1(\R^n)]\times [H^r(\R^n)\cap\, L^1(\R^n)]$,
with $r$ and $s$ defined in each case, the following decay estimates
for the problem {\rm(\ref{ibq})-(\ref{ibqdados})} are true:

\vspace{0,15cm} \noindent (i) If $\,n> 2$ then for $t \geq 0$
$$\begin{array}{ll} ||
v(t)|| \lesssim \left\{||v_1||_{H^r\cap L^1}+ ||v_0||_{H^s\cap
L^1}\right\}
(1+t)^{-\frac{1}{\theta}\left(\frac{n}{4}-\frac{1}{2}\right)},\nonumber
\end{array}$$
with $r=-1$ and $s=0$ in case $\theta =1$, and with
$r=\frac{(1-\theta)(n-2)}{2\theta}-1$ and
$s=\frac{(1-\theta)(n-2)}{2\theta}$ in case $\theta<1$.

\vspace{0,15cm} \noindent (ii) If $\,n\geq 1$ then for $t \geq 0$
$$\begin{array}{ll}
||\nabla v(t)|| + || v_t(t)|| \lesssim \left\{||v_1||_{H^r\cap L^1}+
||v_0||_{H^s\cap L^1}\right\} (1+t)^{-\frac{n}{4\theta}},\nonumber
\end{array}$$
with $r=0$ and $s=1$ in case $\theta =1$, and with
$r=\frac{(1-\theta)n}{2\theta}$ and $s=1+\frac{(1-\theta)n}{2\theta}$
in case $\theta<1$;
$$\begin{array}{ll}
||\Delta v(t)|| + ||\nabla v_t(t)|| \lesssim \left\{||v_1||_{H^r\cap
L^1}+ ||v_0||_{H^s\cap L^1}\right\}
(1+t)^{-\frac{1}{\theta}\left(\frac{n}{4}+\frac{1}{2}\right)},\nonumber
\end{array}$$
with $r=1$ and $s=2$ in case $\theta =1$, and with
$r=1+\frac{(1-\theta)(n+2)}{2\theta}$ and
$s=2+\frac{(1-\theta)(n+2)}{2\theta}$ at case $\theta<1$.
\end{teo}

%%%%%%%%%%%%%%%%%%%%%%%%%%%%%%%%%%%%%%%%%%%%%%%%%
%%%%%%%%%%%%%%%%%%%%%%%%%%%%%%%%%%%%%%%%%%%%%%%%%
%%%%%%%%%%%%%%%% BIBLIOGRAFIA %%%%%%%%%%%%%%%%%%%
%%%%%%%%%%%%%%%%%%%%%%%%%%%%%%%%%%%%%%%%%%%%%%%%%
%%%%%%%%%%%%%%%%%%%%%%%%%%%%%%%%%%%%%%%%%%%%%%%%%

\noindent\textbf{Acknowledgments:} The second author (C. R. L.) was partially supported
by a Research Grant of the Brazilian Government (MCT, Brasil),
Proc. 308868/2015-3.

\end{document}